\documentclass{article}

\usepackage[all]{xy}

\usepackage[text={156mm,244mm},centering]{geometry}
\newlength{\colsep}
\usepackage{multicol}

\usepackage{titlesec}
\titleformat*{\section}{\large\bfseries}
\titleformat*{\subsubsection}{\bfseries}

\usepackage{graphicx}
\usepackage{color}
\usepackage{amsmath}

\definecolor{myurlcolor}{rgb}{0,0,0.5}
\usepackage{hyperref}
\hypersetup{colorlinks,
linkcolor=black,
citecolor=black,
urlcolor=myurlcolor}

\usepackage{bbm}

\usepackage[numbers]{natbib}
\usepackage{natbibspacing}
\renewcommand{\bibfont}{\small}

\usepackage{latexsym}
\usepackage{amssymb}

\newcommand{\slsh}{/\linebreak[0]}
\newcommand{\dblslsh}{//\linebreak[0]}
\newcommand{\dt}{.\linebreak[0]}
\newcommand{\epsln}{\varepsilon}
\newcommand{\integers}{\mathbb{Z}}
\newcommand{\reals}{\mathbb{R}}
\newcommand{\rationals}{\mathbb{Q}}
\newcommand{\complexes}{\mathbb{C}}
\newcommand{\demph}[1]{\textbf{#1}}
\newcommand{\nat}{\mathbb{N}}	
\newcommand{\pr}{\mathrm{pr}}
\newcommand{\of}{\,\raisebox{0.08ex}{\ensuremath{\scriptstyle\circ}}\,}
\newcommand{\sub}{\subseteq}
\newcommand{\One}{\mathbf{1}}
\newcommand{\Two}{\mathbf{2}}
\newenvironment{axiom}[1]%
{\fmlon\begin{trivlist}\item\textbf{#1.}}%
{\end{trivlist}\fmloff}
\newcommand{\axhead}[1]{\minor{#1}}
\newcommand{\minor}[1]{\subsubsection*{#1}}
\newenvironment{defnprop}%
{\begin{trivlist}\item}%
{\end{trivlist}}
\newcommand{\fmlon}{\sffamily}
\newcommand{\fmloff}{\rmfamily}
\newenvironment{dfn}{\medskip\fmlon}{\fmloff\medskip}
\newcommand{\bam}[1]{\begin{center}\it #1\end{center}}
\newcommand{\toby}[1]{\stackrel{#1}{\to}}
\newcommand{\otby}[1]{\stackrel{#1}{\ot}}
\newcommand{\from}{\colon}
\newcommand{\incl}{\hookrightarrow}
\renewcommand{\to}{\longrightarrow}
\newcommand{\ot}{\longleftarrow}
\newcommand{\qsim}{\!\!\sim}
\newcommand{\pset}{\mathcal{P}}
\newcommand{\du}{\sqcup}
\newenvironment{myitemize}%
{\begin{list}%
{$\bullet$}%
{%
\setlength{\leftmargin}{1.5em}%
\setlength{\itemsep}{-.2ex}%
}}%
{\end{list}}
\newcommand{\pn}{}

\title{\Huge\bf Rethinking set theory} 
\author{\Large Tom Leinster}
\date{}

\begin{document}

\sloppy
\maketitle

\setlength{\columnsep}{1\colsep}

\begin{figure*}[b]
\centering
\begin{tabular}{|rl|}
\hline
&\\[-2.2ex]
1       &
Composition of functions is associative and has identities\pn\\
2       &
There is a set with exactly one element\pn      \\
3       &
There is a set with no elements\pn \\
4       &
A function is determined by its effect on elements\pn      \\
5       &
Given sets $X$ and $Y$, one can form their cartesian product $X \times Y$\pn\\
6       &
Given sets $X$ and $Y$, one can form the set of functions from $X$ to $Y$\pn\\
7       &
Given $f\from X \to Y$ and $y \in Y$, one can form the
inverse image $f^{-1}(y)$\pn\\
8       &
The subsets of a set $X$ correspond to the functions from $X$ to $\{0, 1\}$\pn
\\ 
9       &
The natural numbers form a set\pn  \\
10      &
Every surjection has a right inverse\pn   \\[0.3ex]
\hline 
\end{tabular}
\caption{Informal summary of the axioms.  The primitive concepts are set,
function and composition of functions.  Other concepts mentioned
(such as element) are defined in terms of the primitive 
concepts.}
\label{fig:axioms}
\end{figure*}

\begin{multicols}{2}

\begin{trivlist}\item%
\centering%
\parbox{.4\textwidth}{%
\itshape
As mathematicians, we often read a nice new proof of a known theorem, enjoy
the different approach, but continue to derive our internal understanding
from the method we originally learned.  This paper aims to change
drastically the way mathematicians think [\ldots] and teach.%
}
\end{trivlist}
---Sheldon Axler~\cite[Section 10]{AxleDWD}.

\vspace{2em}\noindent
Mathematicians manipulate sets with confidence almost every day of their
working lives.  We do so whenever we work with sets of real or complex
numbers, or with vector spaces, topological spaces, groups, or any of the
many other set-based structures.  These underlying set-theoretic
manipulations are so automatic that we seldom give them a thought, and it
is rare that we make mistakes in what we do with sets.

However, very few mathematicians could accurately quote what are often
referred to as `the' axioms of set theory.  We would not dream of working
with, say, Lie algebras without first learning the axioms.  Yet many of us
will go our whole lives without learning `the' axioms for sets, with no
harm to the accuracy of our work.  This suggests that we all carry around
with us, more or less subconsciously, a reliable body of operating
principles that we use when manipulating sets.  

What if we were to write down some of these principles and adopt
\emph{them} as our axioms for sets?  The message of this article is that
this can be done, in a simple, practical way.  We describe an axiomatization 
due to F.~William Lawvere~\cite{LawvETCS,LawvETCS2}, informally
summarized in Fig.~\ref{fig:axioms}.  The axioms suffice for very nearly
everything mathematicians ever do with sets.  So we can, if we want,
abandon the classical axioms entirely and use these instead.

\minor{Why rethink?}

The traditional axiomatization of sets is known as Zermelo--Fraenkel with
Choice (ZFC).  Great things have been achieved on this axiomatic basis.
However, ZFC has one major flaw: its use of the word `set' conflicts with
how most mathematicians use it.  

The root of the problem is that in the framework of ZFC, the elements of a
set are always sets too.  Thus, given a set $X$, it always makes sense in
ZFC to ask what the elements of the elements of $X$ are.  Now, a typical
set in ordinary mathematics is $\reals$.  But accost a mathematician at
random and ask them `what are the elements of $\pi$?', and they will
probably assume they misheard you, or ask you what you're talking about, or
else tell you that your question makes no sense.  If forced to answer, they
might reply that real numbers have no elements.  But this too is in
conflict with ZFC's usage of `set': if all elements of $\reals$ are sets,
and they all have no elements, then they are all the empty set, from which
it follows that all real numbers are equal.

Some of the actual axioms of ZFC are equally at odds with ordinary
mathematical usage.  For example, one states that every nonempty set $X$
has some element $x$ such that $x \cap X = \emptyset$.  When $X$ is an
ordinary set such as $\reals$, this is a statement that few would recognize
as meaningful: what is $\pi \cap \reals$, after all?

I will anticipate an objection to these criticisms.  The traditional
approach to set theory involves not only ZFC, but also a collection of
methods for encoding mathematical objects of many different types (real
numbers, differential operators, random variables, the Riemann zeta
function, \ldots)\ as sets.  This is similar to the way in which computer
software encodes data of many types (text, sound, images, \ldots)\ as
binary sequences.  In both cases, even the designers would agree that the
encoding methods are somewhat arbitrary.  So, one might object, no one 
is claiming that questions like `what are the elements of $\pi$?' have 
meaningful answers.

However, our understanding that the encoding is not to be taken too
seriously does not alter the bare facts: that in ZFC, it is always valid to
ask of a set `what are the elements of its elements?', and in ordinary
mathematical practice, it is not.  Perhaps it is misleading to use the same
word, `set', for both purposes.

\minor{Three misconceptions}

The axiomatization presented below is Lawvere's Elementary Theory of the
Category of Sets, first proposed half a century ago~\cite{LawvETCS,LawvETCS2}.
Here it is phrased in a way that requires no knowledge of category theory
whatsoever.  

Because of the categorical origins of this axiomatization, three
misconceptions commonly arise.

The first is that the underlying motive is to replace set theory with
category theory.  It is not.  The approach described here is not a rival to
set theory: it \emph{is} set theory.

The second is that this axiomatization demands more mathematical
sophistication than others (such as ZFC).  This is false but
understandable.  Almost all of the work on Lawvere's axioms has taken place
within topos theory: a beautiful and profound subject, but not one easily
accessible to outsiders.  It has always been known that the axioms could be
presented in a completely elementary way, and although some authors have 
emphasized this \cite{LawvETCS,LaRo,MacLMFF,McLaNCB,McLaFOM}, it is not as
widely appreciated as it should be.  This paper aims to make it plain.

The third misconception is that because these axioms for sets come from
category theory, and because the definition of category involves a
collection of objects and a collection of arrows, and because `collection'
might mean something like `set', there is a circularity: in order to
axiomatize sets categorically, we must already know what a set is.  But
although our approach is categorically inspired, it does not depend on
having a general definition of category.  Indeed, our axiomatization
(Section~\ref{sec:axioms}) does not contain a single instance of the word
`category'.

Put another way, circularity is no more a problem here than in ZFC.
Informally, ZFC says `there are some things called sets, there is a binary
relation on sets called membership, and some axioms hold'.  We will say
`there are some things called sets and some things called functions, there
is an operation called composition of functions, and some axioms hold'.  In
neither case are the `things' required to form a set (whatever that would
mean).  In logical terminology, both axiomatizations are simply first-order
theories.

\section{Prelude: elements as functions}
\label{sec:fns-and-elts}

The working mathematician's vocabulary includes terms such as set,
function, element, subset, and equivalence relation.  Any axiomatization of
sets will choose some of these concepts as primitive and derive the others.
The traditional choice is sets and elements.  We use sets and functions.

The formal axiomatization is presented in Section~\ref{sec:axioms}.
However, it will be helpful to consider one aspect in advance: how to
derive the concept of element from the concept of function.

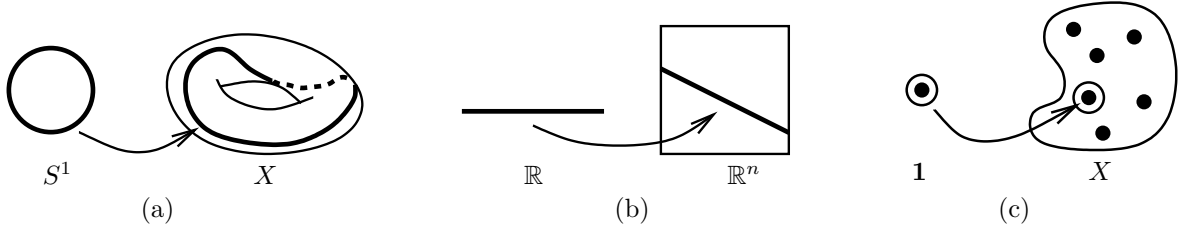
\begin{figure*}
\setlength{\unitlength}{1mm}
\begin{picture}(155,26)(0,-6)
\put(0,0){\includegraphics[width=48\unitlength]{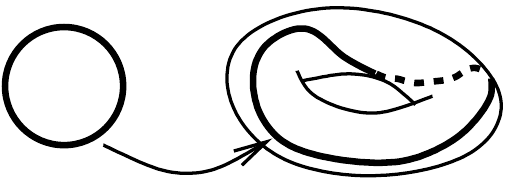}}
\put(5,-3.5){$S^1$}
\put(33,-3.5){$X$}
\put(18,-8){(a)}
\put(60,0){\includegraphics[width=45\unitlength]{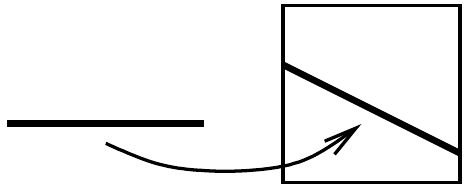}}
\put(69,-3.5){$\reals$}
\put(96,-3.5){$\reals^n$}
\put(81,-8){(b)}
\put(119.5,0.5){\includegraphics[width=37\unitlength]{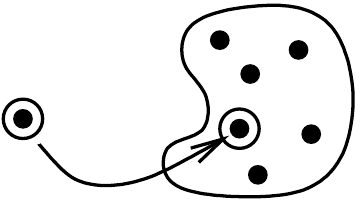}}
\put(120.5,-3){$\One$}
\put(144,-3){$X$}
\put(132,-8){(c)}
\end{picture}
\caption{Mapping out of a basic object ($S^1$, $\reals$, or $\One$) picks out
figures of the appropriate type (loops, lines, or elements).}
\label{fig:elements}
\end{figure*}

Suppose for now that we have found a characterization of one-element sets
without knowing what an element is.  (We do so below.)  Fix a one-element
set $\One = \{ \bullet \}$.  For any set $X$, a function $\One \to X$ is
essentially just an element of $X$, since, after all, such a function $f$
is uniquely determined by the value of $f(\bullet) \in X$
(Fig.~\ref{fig:elements}(c)).  Thus:
\bam{Elements are a special case of functions.}

This is such a trivial observation that one is apt to dismiss it as a mere
formal trick.  On the contrary, similar correspondences occur throughout
mathematics.  For example (Fig.~\ref{fig:elements}):
\begin{myitemize}
\item a loop in a topological space $X$ is a continuous map $S^1 \to X$;
\item a straight line in $\reals^n$ is a distance-preserving map $\reals
  \to \reals^n$;
\item a sequence in a set $X$ is a function $\nat \to X$;
\item a solution $(x, y)$ of the equation $x^2 + y^2 = 1$ in a ring $A$ is
  a homomorphism $\integers[X, Y]/(X^2 + Y^2 - 1) \to A$.  
\end{myitemize}
In each case, the word `is' can be taken either as a definition or as an
assertion of a canonical, one-to-one correspondence.  In the first, we map
out of the circle, which is a `free-standing' loop; in the second, $\reals$
is a free-standing line; in the third, the elements $0, 1, 2, \ldots$\ of
$\nat$ form a free-standing sequence; in the last, the pair $(X, Y)$ of
elements of $\integers[X, Y]/(X^2 + Y^2 - 1)$ is the free-standing solution
$(x, y)$ of $x^2 + y^2 = 1$.  Similarly, in our trivial situation, the set
$\One$ is a free-standing element, and an element of a set $X$ is just a
map $\One \to X$.

We could be fussy and write $\bar{x}$, say, for the function $\One \to X$
with value $x \in X$.  But we will write $\bar{x}$ as just $x$, blurring
the distinction.  In fact, we will later \emph{define} an element of $X$ to
be a function $\One \to X$.

This will make some readers uncomfortable.  There is, you will agree, a
canonical one-to-one correspondence between elements of $X$ and functions
$\One \to X$, but perhaps you draw the line at saying that an element of
$X$ literally \emph{is} a function $\One \to X$.  If so, this is not a
deal-breaker.  We could adapt the axiomatization
in Section~\ref{sec:axioms} 
by adding `element' to the list of primitive
concepts.  Then, however, we would need to complicate it further by adding
clauses to guarantee that (among other things) there is a one-to-one
correspondence between elements of $X$ and functions $\One \to X$, for any
set $X$.  It can be done, but we choose the more economical route.

We have seen that elements are a special case of functions.  There is
another fundamental way in which functions and elements interact: given a
function $f\from X \to Y$ and an element $x \in X$, we can evaluate $f$ at
$x$ to obtain a new element, $f(x) \in Y$.  Viewing elements as functions
out of $\One$, this element $f(x)$ is nothing but the composite of $f$ with
$x$:
\[
f(x) = f \of x.
\qquad\qquad
\begin{array}{c}
\xymatrix{
\One \ar[r]^x \ar[rd]_{f(x)}    &X \ar[d]^f     \\
                                &Y
}
\end{array}
\]
Hence:
\bam{Evaluation is a special case of composition.}

\section{The axioms} 
\label{sec:axioms}

Here we state our ten axioms on sets and functions, in entirely elementary
terms.

The formal axiomatization is in a \fmlon different typeface\fmloff, to
distinguish it from the accompanying commentary.  Some diagrams appear, but
they are not part of the formal statement.

First we state the data to which our axioms will apply:
\fmlon
\begin{myitemize}
\item Some things called \demph{sets};

\item for each set $X$ and set $Y$, some things called \demph{functions
from $X$ to $Y$}, with functions $f$ from $X$ to $Y$ written as $f\from
X \to Y$ or $X \toby{f} Y$;

\item for each set $X$, set $Y$ and set $Z$, an operation assigning to each
$f\from X \to Y$ and $g\from Y \to Z$ a function $g\of f\from X \to Z$;

\item for each set $X$, a function $1_X\from X \to X$.
\end{myitemize}
\fmloff

This last item can be included in the list or not, according to taste.  See
the comments after the first axiom, which now follows.

\axhead{Associativity and identity laws} 
\begin{axiom}{1}
For all sets $W, X, Y, Z$ and functions 
\[
W \toby{f} X \toby{g} Y \toby{h} Z,
\]
we have $h \of (g\of f) = (h \of g) \of f$.  For all sets $X, Y$ and functions
$f\from X \to Y$, we have $f \of 1_X = f = 1_Y \of f$.
\end{axiom}

If we wish to omit the identity functions from the list of primitive
concepts, we must replace the second half of Axiom~1 by the statement
that for all sets $X$, there exists a function $1_X \from X \to X$ such
that $g \of 1_X = g$ for all $g\from X \to Y$ and $1_X \of f = f$ for all
$f\from W \to X$.  These conditions characterize $1_X$ uniquely.

\axhead{One-element set}
We would like to say `there exists a one-element set', but for the moment
we lack the expressive power to say `element'.  However, any one-element
set $T$ should have the property that for each set $X$, there is precisely
one function $X \to T$.  Moreover, \emph{only} one-element sets should have
this property.  This motivates the following definition and axiom.

\begin{dfn}
A set $T$ is \demph{terminal} if for every set $X$, there is a unique
function $X \to T$.
\end{dfn}

\begin{axiom}{2}
There exists a terminal set.
\end{axiom}

It follows quickly from the definitions that if $T$ and $T'$ are terminal
sets then there is a unique isomorphism from $T$ to $T'$.  (A function $f
\from A \to B$ is an \demph{isomorphism} if there is a function $f' \from B
\to A$ such that $f' \of f = 1_A$ and $f \of f' = 1_B$.)  In other words,
terminal sets are unique up to unique isomorphism.  It is therefore
harmless to fix a terminal set $\One$ once and for all.  Readers worried by
this are referred to the last few paragraphs of this section.

\begin{dfn}
Given a set $X$, we write $x \in X$ to mean $x\from \One \to X$, and call
$x$ an \demph{element} of $X$.  Given $x \in X$ and a function $f\from X
\to Y$, we write $f(x)$ for the element $f \of x\from \One \to Y$ of
$Y$. 
\end{dfn}

\axhead{Empty set}

\begin{axiom}{3}
There exists a set with no elements.
\end{axiom}

\axhead{Functions and elements}
A function from $X$ to $Y$ should be nothing more than a way of turning
elements of $X$ into elements of $Y$.  

\begin{axiom}{4}
Let $X$ and $Y$ be sets and $f, g\from X \to Y$ functions.  Suppose that
$f(x) = g(x)$ for all $x \in X$.  Then $f = g$.
\end{axiom}

Axioms~1, 2 and~4 imply that a set is terminal if and only if it has
exactly one element.  This justifies the usage of `one-element set' as a
synonym for `terminal set'.

\axhead{Cartesian products}
We want to be able to form cartesian products of sets.  An element of $X$
together with an element of $Y$ should uniquely determine an element of $X
\times Y$.  More generally, for any set $I$, a function $f_1 \from I \to X$
together with a function $f_2 \from I \to Y$ should uniquely determine a
function $f \from I \to X \times Y$, given by $f(t) = (f_1(t),
f_2(t))$.  (To see that this really is `more generally', take $I =
\One$.)  We can recover $f_1$ from $f$ by composing with the projection
$p_1 \from X \times Y \to X$, and similarly $f_2$, as in the following
definition. 

\begin{dfn}
Let $X$ and $Y$ be sets.  A \demph{product} of $X$ and $Y$ is a set $P$
together with functions $X \otby{p_1} P \toby{p_2} Y$, with the following
property:
\begin{defnprop}
for all sets $I$
and functions $X \otby{f_1} I \toby{f_2} Y$,\\
there is a unique function $(f_1, f_2)\from I \to P$\\
such that
$p_1 \of (f_1, f_2) = f_1$ and $p_2 \of (f_1, f_2) = f_2$.
\end{defnprop}
\end{dfn}
\[
\begin{xy}
(0,20)*+{I}="test";
(0,8)*+{P}="prod";
(-20,0)*+{X}="first";
(20,0)*+{Y}="second";
{\ar_{f_1} "test"; "first"};
{\ar^{f_2} "test"; "second"};
{\ar@{.>}|{(f_1, f_2)} "test"; "prod"};
{\ar^{p_1} "prod"; "first"};
{\ar_{p_2} "prod"; "second"};
\end{xy}
\]
\begin{axiom}{5}
Every pair of sets has a product.
\end{axiom}

Strictly speaking, a product consists of not only the set $P$ but also the
projections $p_1$ and $p_2$.  Any two products of $X$ and $Y$ are uniquely
isomorphic: that is, given products $(P, p_1, p_2)$ and $(P', p'_1, p'_2)$,
there is a unique isomorphism $i\from P \to P'$ such that $p'_1 \of i =
p_1$ and $p'_2 \of i = p_2$.  As in the case of terminal sets, this makes
it harmless to choose once and for all a preferred product $(X \times Y,
\pr^{X, Y}_1, \pr^{X, Y}_2)$ for each pair $X$, $Y$ of sets.  Again, this
convention is justified at the end of the section.

\axhead{Sets of functions}
In everyday mathematics, we can form the set $Y^X$ of functions from one
set $X$ to another set $Y$.  For any set $I$, the functions $q \from I
\times X \to Y$ correspond one-to-one with the functions $\bar{q}\from I
\to Y^X$, simply by changing the punctuation:
\begin{equation}        \label{eq:curry}
q(t, x)
=
(\bar{q}(t))(x)
\end{equation}
($t \in I$, $x \in X$).  For example, when $I = \One$, this reduces to the
statement that the functions $X \to Y$ correspond to the elements of $Y^X$.

In~\eqref{eq:curry}, we are implicitly using the evaluation map
\[
\begin{array}{cccc}
\epsln\from     &Y^X \times X   &\to            &Y      \\
                &(f, x)         &\longmapsto    &f(x).
\end{array}
\]
Then~\eqref{eq:curry} becomes the equation $q(t, x) = \epsln (\bar{q}(t),
x)$, as in the following definition.

\begin{dfn}
Let $X$ and $Y$ be sets.  A \demph{function set} from $X$ to $Y$ is a
set $F$ together with a function $\epsln\from F \times X \to Y$, with the
following property: 
\begin{defnprop}
for all sets $I$ and functions $q\from I \times X \to Y$,\\
there is a unique function $\bar{q}\from I \to F$\\
such that 
$q(t, x) = \epsln (\bar{q}(t), x)$ for all $t \in I$, $x \in X$.
\end{defnprop}
\end{dfn}
\[
\xymatrix{
I \times X \ar@{.>}[d]_{\bar{q} \times 1_X} \ar[rd]^q   &       \\
F \times X \ar[r]_\epsln                                &Y
}
\]
\begin{axiom}{6}
For all sets $X$ and $Y$, there exists a function set from $X$ to $Y$.
\end{axiom}

\axhead{Inverse images}
Ordinarily, given a function $f\from X \to Y$ and an element $y$ of $Y$, we
can form the inverse image or fibre $f^{-1}(y)$.  The inclusion function $j
\from f^{-1}(y) \incl X$ has the property that $f\of j$ has constant value
$y$.  Moreover, whenever $q \from I \to X$ is a function such that $f \of
q$ has constant value $y$, the image of $q$ must lie within $f^{-1}(y)$;
that is, $q = j \of \bar{q}$ for some $\bar{q}\from I \to f^{-1}(y)$
(necessarily unique).

\begin{dfn}
Let $f\from X \to Y$ be a function and $y \in Y$.  An \demph{inverse image}
of $y$ under $f$ is a set $A$ together with a function $j\from A \to X$, such
that $f(j(a)) = y$ for all $a \in A$ and the following property holds:
\begin{defnprop}
for all sets $I$ and functions $q \from I \to X$ such that\\
\mbox{}\quad $f(q(t)) = y$ for all $t \in I$,\\
there is a unique function $\bar{q} \from I \to A$\\ 
such that $q = j \of \bar{q}$.  
\end{defnprop}
\end{dfn}
\[
\begin{xy}
(0,12)*+{A}="inv";
(12,12)*+{\One}="one";
(0,0)*+{X}="dom";
(12,0)*+{Y}="cod";
(-13,18)*+{I}="test";
{\ar^j "inv"; "dom"};
{\ar "inv"; "one"};
{\ar_f "dom"; "cod"};
{\ar^y "one"; "cod"};
{\ar@/_/_{q} "test"; "dom"};
{\ar@/^/ "test"; "one"};
{\ar@{.>}|{\bar{q}} "test"; "inv"};
\end{xy}
\]

\begin{axiom}{7}
For every function $f\from X \to Y$ and element $y \in Y$, there exists an
inverse image of $y$ under $f$.
\end{axiom}

Inverse images are essentially unique: if $j\from A \to X$ and
$j'\from A' \to X$ are both inverse images of $y$ under $f$, there is a
unique isomorphism $i \from A \to A'$ such that $j' \of i = j$.

\axhead{Characteristic functions}
Sometimes we want to define a function on a case-by-case basis.  For
example, we might want to define $h\from \reals \to \reals$ by $h(x) = x
\sin(1/x)$ if $x \neq 0$ and $h(0) = 0$.  A simple instance is the
definition of characteristic function.

Fix a two element-set $\Two = \{t, f\}$ (for `true' and `false').  The
characteristic function of a subset $A \sub X$ is the function $\chi_A
\from X \to \Two$ defined by $\chi_A(x) = t$ if $x \in A$ and $\chi_A(x)
= f$ otherwise.  It is the unique function $\chi\from X \to \Two$ such
that $\chi^{-1}(t) = A$.  

This is how characteristic functions work ordinarily.  To ensure that they
work in the same way in our set theory, we now demand that there exist a
set $\Two$ and an element $t \in \Two$ with the property just described:
whenever $X$ is a set and $A \sub X$, there is a unique function $\chi\from
X \to \Two$ such that $\chi^{-1}(t) = A$.

Since we do not yet have a definition of subset, we phrase the axiom in
terms of injections instead.  This works because every subset inclusion $A
\incl X$ is injective, and, up to isomorphism, every injection arises in
this way.

\begin{dfn}
An \demph{injection} is a function $j \from A \to X$ such that $j(a) =
j(a') \implies a = a'$ for $a, a' \in A$.

A \demph{subset classifier} is a set $\Two$ together with an element $t \in
\Two$, with the following property:
\begin{defnprop}
for all sets $A, X$ and injections $j\from A \to X$,\\ 
there is a unique function $\chi\from X \to \Two$ such that\\ 
$j\from A \to X$ is an inverse image of $t$ under $\chi$.
\end{defnprop}
\end{dfn}
\[
\xymatrix{
A \ar[r] \ar[d]_j       &\One \ar[d]^t  \\
X \ar@{.>}[r]_\chi      &\Two
}
\]
\begin{axiom}{8}
There exists a subset classifier.
\end{axiom}

The notation $\Two$ is merely suggestive.  There is nothing in the
definition saying that $\Two$ must have two elements, but, nontrivially,
our ten axioms do in fact imply this.

\axhead{Natural numbers}
In ordinary mathematics, sequences can be defined recursively: given a set
$X$, an element $a \in X$, and a function $r \from X \to X$, there is a
unique sequence $(x_n)_{n = 0}^\infty$ in $X$ such that
\[
x_0 = a
\text{ and }
x_{n + 1} = r(x_n) 
\text{ for all } n \in \nat.
\]
A sequence in $X$ is nothing but a function $\nat \to X$, so the previous
sentence is really a statement about the set $\nat$.  It also refers to two
pieces of structure on $\nat$: the element $0$ and the function $s\from
\nat \to \nat$ given by $s(n) = n + 1$.

\begin{dfn}
A \demph{natural number system} is a set $N$ together with an
element $0 \in N$ and a function $s\from N \to N$, with the
following property:
\begin{defnprop}
whenever $X$ is a set, $a \in X$, and $r\from X \to X$,\\
there is a unique function $x\from N \to X$ such that\\
$x(0) = a$
and 
$x(s(n)) = r(x(n))$ for all $n \in N$.
\end{defnprop}
\end{dfn}
\[
\xymatrix{
\One \ar[r]^0 \ar[d]_{1_{\One}} &N \ar[r]^s \ar@{.>}[d]^x
                                                &N \ar@{.>}[d]^x        \\
\One \ar[r]_a                   &X \ar[r]_r     &X
}
\]
\begin{axiom}{9}
There exists a natural number system.
\end{axiom}

Natural number systems are essentially unique, in the usual sense that
between any two of them there is a unique structure-preserving isomorphism.
This justifies speaking of \emph{the} natural numbers $\nat$, as we
invariably do.

\axhead{Choice}
A function with a right inverse is certainly surjective.  The axiom of
choice states the converse.

\begin{dfn}
A \demph{surjection} is a function $s \from X \to Y$ such that for all $y
\in Y$, there exists $x \in X$ with $s(x) = y$.

A \demph{right inverse} of a function $s\from X \to Y$ is a function
$i\from Y \to X$ such that $s \of i = 1_Y$.
\end{dfn}

\begin{axiom}{10}
Every surjection has a right inverse.
\end{axiom}

A right inverse of a surjection $s\from X \to Y$ is a choice, for each $y \in
Y$, of an element of the nonempty set $s^{-1}(y)$.  

This concludes the axiomatization.

\minor{The meaning of `the'}

It remains to reassure any readers concerned by the liberty taken in
Axioms~2 and~5, where we chose once and for all a terminal set and a
cartesian product for each pair of sets.

This type of liberty is very common in mathematical practice.  We speak of
\emph{the} trivial group, \emph{the} 2-sphere, \emph{the} direct sum of two
vector spaces, etc., even though we can conceive of many trivial groups or
2-spheres or direct sums, all isomorphic but not equal.  Anyone asking `but
\emph{which} trivial group?'\ is likely to be met with a hard stare, for
good reason: no meaningful statement about groups depends on what the
element of the trivial group happens to be called.

However, we should be able to state the axioms with scrupulous rigour, and
we can.  One way to do so is not to single out a particular terminal set or
particular products, but instead to adopt some circumlocutions: for
example, replacing the phrase `for all elements $x \in X$' by `for all
terminal sets $T$ and functions $x\from T \to X$'.

More satisfactory, though, is to extend the list of primitive concepts.  To
the existing list (sets, functions, composition and identities) we add:
\fmlon
\begin{myitemize}

\item a distinguished set, $\One$;

\item an operation assigning to each pair of sets $X, Y$ a set $X
\times Y$ and functions
\begin{equation}
\label{eq:prod-pjns}
\xymatrix@1@C+2ex{
X       &
\ar[l]_-{\pr^{X, Y}_1} X \times Y \ar[r]^-{\pr^{X, Y}_2}  &
Y.
}
\end{equation}

\end{myitemize}
\fmloff
Axiom~2 is replaced by the statement that $\One$ is terminal, and
Axiom~5 by the statement that for all sets $X$ and $Y$, the set $X \times
Y$ together with the functions~\eqref{eq:prod-pjns} is a product of $X$ and
$Y$.  

This approach has the virtue of reflecting ordinary mathematical usage.  We
usually speak as if taking the product of two sets (or spaces, groups,
etc.)\ were a procedure with a definite output: \emph{the} product, not
\emph{a} product.  But since products are in any case determined uniquely
up to unique isomorphism, whether or not we nominate one as special makes
no significant difference.

\section{Discussion} 
\label{sec:rmks}

The ten axioms are familiar in their intuitive content, but less so
as an axiomatic system.  Here we discuss the implications of using them as
such. 

\minor{Building on the axioms}

Any axiomatization of anything is followed by a period of lemma-proving.
The present axioms are no exception.  Here is a very brief sketch of the
development.

It is convenient formally to define a \demph{subset} of a set $X$ as a
function $X \to \Two$, but we constantly use the correspondence between
functions $X \to \Two$ and injections into $X$, provided by Axiom~8.  Two
injections $j$, $j'$ into $X$ correspond to the same subset of $X$ if and
only if they have the same image (that is, there exists an isomorphism $i$
such that $j' = j \of i$).

The main task is to build the everyday equipment used for manipulating
sets.  For example, given a function $f\from X \to Y$, we construct the
image under $f$ of a subset of $X$ and the inverse image of a subset of
$Y$.  An equivalence relation $\sim$ on a set $X$ is defined to be a subset
of $X \times X$ with the customary properties, and the axioms allow us to
construct the quotient set $X/\qsim$.  Some constructions are tricky: for
instance, the axioms imply that any two sets $X$ and $Y$ have a disjoint
union $X \du Y$, but this is by no means obvious.

We then define the usual number systems.  Addition, multiplication and
powers of natural numbers are defined directly using Axiom~9.  From $\nat$
we successively construct $\integers$, $\rationals$, $\reals$ and
$\complexes$, in the standard way.  For example, $\integers = (\nat \times
\nat)/\qsim$, where $\sim$ is the equivalence relation on $\nat \times
\nat$ given by $(m, n) \sim (m', n')$ if and only if $m + n' = m' + n$.  As
this illustrates, past a certain point, the development is literally
identical to that for other axiomatizations of sets.

\minor{How strong are the axioms?}

Most mathematicians will never use more properties of sets than those
guaranteed by the ten axioms.  For example, McLarty~\cite{McLaFOA} argues
that no more is needed anywhere in the canons of the Grothendieck school of
algebraic geometry, the multi-volume works \emph{\'El\'ements de
G\'eom\'etrie Alg\'ebrique} (EGA) and \emph{S\'eminaire de G\'eom\'etrie
Alg\'ebrique} (SGA).

To get a sense of the reach of the axioms, let us consider infinite
cartesian products.  Let $I$ be a (possibly infinite) set and $(X_i)_{i \in
  I}$ a family of sets.  Can we form the product $\prod_{i \in I} X_i$?
It depends on what is meant by `family'.  We could define an
$I$-indexed family to be a set $X$ together with a function $p\from X \to
I$, viewing the fibre $p^{-1}(i)$ as the $i$th member $X_i$.  In that case,
$\prod X_i$ can be constructed as a subset of $X^I$.  Specifically, $p$
induces a function $p^I \from X^I \to I^I$, and $\prod X_i$ is the inverse
image under $p^I$ of the element of $I^I$ corresponding to $1_I$.

However, we could interpret `$I$-indexed family' differently: as an
algorithm or formula that assigns to each $i \in I$ a set $X_i$.  It is not
obvious that we can then form the disjoint union $X = \coprod_{i \in I}
X_i$, which is what would be necessary in order to obtain a family in the
previous sense.  In fact, writing $\pset(S) = \Two^S$ for the power set of
a set $S$, the ten axioms do \emph{not} guarantee the existence of the
disjoint union
\begin{equation}        \label{eq:aleph_omega}
\nat \du \pset(\nat) \du \pset(\pset(\nat)) 
\du \cdots 
\end{equation}
unless they are inconsistent (\cite{MathSML}, Section~9).

If we wish to change this, we can add an eleventh axiom (or properly, axiom
scheme), called `replacement' and informally stated as follows.  Suppose we
have a set $I$ and a first-order formula that for each $i \in I$ specifies
a set $X_i$ up to isomorphism.  Then we require that there exist a set $X$
and a function $p \from X \to I$ such that $p^{-1}(i)$ is isomorphic to
$X_i$ for each $i \in I$.  (See Section~8 of~\cite{McLaECS} for a formal
statement.)  This guarantees the existence of sets such
as~\eqref{eq:aleph_omega}.

The relationship between our axioms and ZFC is well understood.  The ten
axioms are weaker than ZFC; but when the eleventh is added, the two
theories have equal strength and are `bi-interpretable' (the same theorems
hold).  Moreover, it is known to which fragment of ZFC the ten axioms
correspond: `Zermelo with bounded comprehension and choice'.  The details
of this relationship were mostly worked out in the early 1970s
\cite{Cole,MitcBTT,Osiu}.  Good modern accounts are in Section~VI.10
of~\cite{MaMo} and Chapter~22 of~\cite{McLaECET}.

\minor{A broader view}

Our ten axioms are a standard rephrasing of Lawvere's Elementary Theory of
the Category of Sets (ETCS), published in 1964.  It was some years before
ETCS found its natural home, and that was with the advent of topos theory.

The notion of topos was invented by Grothendieck for reasons that had
nothing to do with set theory.  For Grothendieck, a topos was a generalized
topological space.  Formally, a topos is a category with
certain properties, and a topological space $X$ is associated with the
topos whose objects are the sheaves of sets on $X$.  

Lawvere and Tierney swiftly realized that, after a slight loosening of
Grothendieck's definition, the ETCS axioms could be restated neatly in
topos-theoretic terms~\cite{TierSTC,TierAST}.  Indeed, ETCS says exactly
that sets and functions form a topos of a special sort: a `well-pointed
topos with natural numbers object and choice'.  So a topos is not only a
generalized space; it is also a generalized universe of sets.

An attractive feature of ETCS is that each of the axioms is meaningful in a
broader context than set theory.  For example, Axiom~1 states that sets and
functions form a category.  The job of the remaining axioms is to
distinguish sets from other structures that form categories.  Axioms~2
and~5 state that the category of sets has finite products.  This important
property is shared by (for example) the categories of topological spaces
and smooth manifolds, which is exactly what makes it possible to define
`topological group' and `Lie group'.  But for one detail, Axioms 1, 2, 5,
6, 7 and~8 state that sets and functions form a topos.

Skipping to Axiom~10, the axiom of choice as formulated there highlights a
special feature of sets.  In most other categories of sets-with-structure,
it fails, and its failure is a point of interest.  For instance, not every
continuous surjection between topological spaces has a continuous right
inverse, a typical example being the nonexistence of a continuous square
root defined on the complex plane.

\minor{What kind of set theory should we teach?}

As Fig.~\ref{fig:axioms} indicates, we already teach a diluted form of the
ten axioms, even in introductory courses.  For example, we certainly tell
our students that an element of $X \times Y$ is an element of $X$ together
with an element of $Y$, and we routinely write a function $f$ taking values
in $\reals^2$ as $(f_1, f_2)$, although we are less likely to state
explicitly that given functions $f_1\from I \to X$ and $f_2\from I \to Y$,
there is a unique function $f\from I \to X \times Y$ with $f_1$ and $f_2$
as components.

When it comes to teaching \emph{axiomatic} set theory, the approach
outlined here has advantages and disadvantages.  The big advantage is that
such a course is of far wider benefit than one using the traditional
axioms.  It directly addresses a difficulty experienced by many students:
the concept of function (and worse, function space).  It also introduces in
an elementary setting the idea of universal property.  This is probably the
hardest aspect of the axioms for a learner, but since universal properties
are important in so many branches of advanced mathematics, the benefits are
potentially far-reaching.

The disadvantages are perhaps only temporary.  There is at present a lack
of teaching materials (the book~\cite{LaRo} being the main exception).  For
example, the axioms imply that any two sets have a disjoint union, and
most books on topos theory contain an elegant and sophisticated proof of a
generalization of this fact, but to my knowledge, there is only one place where
a purely elementary proof can be found~\cite{TrimETCS3}.  A second
disadvantage is that any student planning a career in set theory will need
to learn ZFC anyway, since almost all research-level set theory is done
with the iterated-membership conception of set.  (That is the current
reality, which is not to say that set theory \emph{must} be done this way.)

\minor{Reactions to an earthquake}

Perhaps you will wake up tomorrow, check your email, and find an
announcement that ZFC is inconsistent.  Apparently, someone has taken the
ZFC axioms, performed a long string of logical deductions, and arrived at a
contradiction.  The work has been checked and re-checked.  There is no
longer any doubt.

How would you react?  In particular, how would you feel about the
implications for your own work?  All your theorems would still be true
under ZFC, but so too would their negations.  Would you conclude that your
life's work had been destroyed?

I believe that most of us would be interested but not deeply troubled.
We would go on believing that our theorems were true in a sense that their
negations were not.  We are unlikely to feel threatened by the
inconsistency of axioms to which we never referred anyway.

In contrast, the ten axioms above are such core mathematical principles
that an inconsistency in them would be devastating.  If we cannot safely
assume that composition of functions is associative, or that repeatedly
applying a function $f\from X \to X$ to an element $a \in X$ produces a
sequence $(f^n(a))$, we are really in trouble.

As the weaker system, the ten axioms are less likely to be inconsistent
than ZFC.  But the question of strength is peripheral to this article (and
in any case, if one wants a system of equal strength to ZFC, all one needs
to do is add the aforementioned eleventh axiom).  The real message is this:
simply by writing down a few mundane, uncontroversial statements about sets
and functions, we arrive at an axiomatization that reflects how sets are
used in everyday mathematics.

\paragraph*{Acknowledgements}
I thank Fran\c{c}ois Dorais, Colin McLarty, Todd Trimble and the patrons of
the $n$-Category Caf\'e.  This work was partially supported by an EPSRC
Advanced Research Fellowship.

\bigskip
\noindent
School of Mathematics,
University of Edinburgh,\\
Edinburgh EH9 3JZ,
United Kingdom.\\
Tom.Leinster\mbox{}@\mbox{}ed\dt ac\dt uk

\end{multicols}


\begin{thebibliography}{10}

\bibitem{AxleDWD}
S.~Axler.
\newblock Down with determinants!
\newblock {\em American Mathematical Monthly}, 102:139--154, 1995.

\bibitem{Cole}
J.~C. Cole.
\newblock Categories of sets and models of set theory.
\newblock In J.~Bell and A.~Slomson, editors, {\em Proceedings of the Bertrand
  Russell Memorial Logic Conference, Uldum 1971}, pages 351--399. 1973.

\bibitem{LawvETCS}
F.~W. Lawvere.
\newblock An elementary theory of the category of sets.
\newblock {\em Proceedings of the National Academy of Sciences of the U.S.A.},
  52:1506--1511, 1964.

\bibitem{LawvETCS2}
F.~W. Lawvere.
\newblock An elementary theory of the category of sets (long version) with
  commentary.
\newblock {\em Reprints in Theory and Applications of Categories}, 12:1--35,
  2005.

\bibitem{LaRo}
F.~W. Lawvere and R.~Rosebrugh.
\newblock {\em Sets for Mathematics}.
\newblock Cambridge University Press, Cambridge, 2003.

\bibitem{MacLMFF}
S.~Mac~Lane.
\newblock {\em Mathematics: Form and Function}.
\newblock Springer, New York, 1986.

\bibitem{MaMo}
S.~Mac~Lane and I.~Moerdijk.
\newblock {\em Sheaves in Geometry and Logic}.
\newblock Springer, New York, 1994.

\bibitem{MathSML}
A.~Mathias.
\newblock The strength of {M}ac {L}ane set theory.
\newblock {\em Annals of Pure and Applied Logic}, 110:107--234, 2001.

\bibitem{McLaECET}
C.~Mc{L}arty.
\newblock {\em Elementary Categories, Elementary Toposes}.
\newblock Oxford University Press, 1992.

\bibitem{McLaNCB}
C.~Mc{L}arty.
\newblock Numbers can be just what they have to.
\newblock {\em No{\^{u}}s}, 27:487--98, 1993.

\bibitem{McLaFOM}
C.~Mc{L}arty.
\newblock Challenge axioms, final draft.
\newblock Email to Foundations of Mathematics mailing list, 6 February 1998,
  archived at \href{http://www.cs.nyu.edu/pipermail/fom}{http:\dblslsh www\dt
  cs\dt nyu\dt edu\slsh pipermail\slsh fom}, 1998.

\bibitem{McLaECS}
C.~Mc{L}arty.
\newblock Exploring categorical structuralism.
\newblock {\em Philosophia Mathematica}, 12:37--53, 2004.

\bibitem{McLaFOA}
C.~Mc{L}arty.
\newblock A finite order arithmetic foundation for cohomology.
\newblock \href{http://arxiv.org/abs/1102.1773}{arXiv:1102.1773}, 2011.

\bibitem{MitcBTT}
W.~Mitchell.
\newblock Boolean topoi and the theory of sets.
\newblock {\em Journal of Pure and Applied Algebra}, 2:261--274, 1972.

\bibitem{Osiu}
G.~Osius.
\newblock Categorical set theory: a characterization of the category of sets.
\newblock {\em Journal of Pure and Applied Algebra}, 4:79--119, 1974.

\bibitem{TierSTC}
M.~Tierney.
\newblock Sheaf theory and the continuum hypothesis.
\newblock In F.~W. Lawvere, editor, {\em Toposes, Algebraic Geometry and
  Logic}, volume 274 of {\em Lecture Notes in Mathematics}, pages 13--42.
  Springer, 1972.

\bibitem{TierAST}
M.~Tierney.
\newblock Axiomatic sheaf theory: some constructions and applications.
\newblock In P.~Salmon, editor, {\em Proceedings of CIME Conference on
  Categories and Commutative Algebra, Varenna, 1971}, pages 249--326. Edizione
  Cremonese, 1973.

\bibitem{TrimETCS3}
T.~Trimble.
\newblock {ETCS}: building joins and coproducts.
\newblock \href{http://ncatlab.org/nlab/show/Trimble+on+ETCS+III}{http:\dblslsh
  ncatlab.org\slsh nlab\slsh show\slsh
  Trimble+\linebreak[0]on+\linebreak[0]ETCS+III}, 2008.

\end{thebibliography}
\end{document}